\documentclass[11pt,leqno]{article}

\usepackage{amsfonts,latexsym, amsmath, amssymb, amsthm, hyperref}

\usepackage{gfsartemisia}
\usepackage[T1]{fontenc}

\usepackage{color}
\usepackage{fancyhdr} 
\usepackage{lastpage}

\newtheorem{theorem}{Theorem}[section]

\newtheorem{proposition}[theorem]{Proposition}
\newtheorem{corollary}[theorem]{Corollary}
\newtheorem{remark}[theorem]{Remark}
\newtheorem{remarks}[theorem]{Remarks}

\theoremstyle{definition}

\theoremstyle{remark}
\newtheorem*{note*}{Note}

\makeatletter

\numberwithin{equation}{section}

\makeatother

\newcommand{\prend}{$\hfill \quad \Box$}

\newcommand\blfootnote[1]{%
  \begingroup
  \renewcommand\thefootnote{}\footnote{#1}%
  \addtocounter{footnote}{-1}%
  \endgroup
}

\begin{document}

\small

\title{A Gaussian small deviation inequality for convex functions}

\author{Grigoris Paouris\thanks{Supported by the NSF
CAREER-1151711 grant;} \, and \, Petros Valettas\thanks{Supported in part by the NSF grant DMS-1612936.}}

%\date{}

\maketitle

\begin{abstract}\footnotesize
Let $Z$ be an $n$-dimensional Gaussian vector and let $f:\mathbb R^n \to \mathbb R$ be a convex function. We prove that:
\begin{align*}
\mathbb P \left( f(Z) \leq \mathbb Ef(Z) -t\sqrt{ {\rm Var}f(Z) } \right)\leq \exp(-ct^2),
\end{align*} for all $t>1$ where $c>0$ is an absolute constant. 
 As an application we derive variance-sensitive small ball probabilities for Gaussian processes.
\end{abstract}


%%%%%%%%%%%%
\section{Introduction}
%%%%%%%%%%%%

The purpose of this note is to establish a sharp distributional inequality for convex functions on Gauss' space 
$(\mathbb R^n, \|\cdot\|_2, \gamma_n)$. 
Our goal and motivation stems from the attempt to strengthen the classical Gaussian concentration for special cases that are of interest in 
high-dimensional
geometry. The Gaussian concentration phenomenon (see \cite{Bog} and \cite{Led}) states 
that for any $L$-Lipschitz map $f: \mathbb R^n \to \mathbb R$ one has
\begin{align} \label{eq:Gauss-conc}
\mathbb P\left( \big| f(Z) - M \big| >t \right) \leq \exp \left( -\tfrac{1}{2}t^2/L^2 \right),
\end{align} for all $t>0$, where $Z$ is $n$-dimensional standard Gaussian random vector and $M$ is a median for $f(Z)$. 
The above inequality follows from the solution to the isoperimetric problem in Gauss' space, which was proved independently
by Borell in \cite{Bor-iso} and Sudakov and Tsirel'son in \cite{ST} and can be described by the following inequality:
\begin{align} \label{eq:iso-ineq}
\gamma_n(A + tB_2^n) \geq \Phi( \Phi^{-1}(\gamma_n(A))+t), \quad t>0,
\end{align} for all Borel sets $A\subseteq \mathbb R^n$, where $\Phi$ is the cumulative distribution
function of a Gaussian random variable. Applying \eqref{eq:iso-ineq} for 
$A=\{f\leq M\}$, where $M$ is a median of $f$ and by taking into account that $A+tB_2^n \subseteq \{f\leq M +tL \}$ we obtain:
\begin{align*}
\gamma_n(f\leq M +tL) \geq \Phi(t) \quad  \Longrightarrow \quad \gamma_n( f > M + tL ) \leq 1-\Phi(t).
\end{align*} Finally, standard estimates for the function $\Phi$, such as \eqref{eq:ineq-phi}, yield the result (we work similarly for the deviation
below the median). In turn this implies bounds 
on the variance ${\rm Var}f(Z)$ for any Lipschitz map $f$ in terms of the Lipschitz constant $L$ (alternatively we may employ 
the Gaussian Poincar\'e inequality \cite{Chen}):
\begin{align} \label{eq:var-ineq}
{\rm Var}[f(Z)] \leq  L^2.
\end{align} The above inequalities are sharp for linear functionals. However, one can easily construct examples of convex functions 
(see Section 2) for which the above estimates are far from being optimal. 
On the other hand the observation in \cite[Corollary 3.2]{LedTal} that for $f$ being a norm, one has the 
stochastic dominance $\mathbb P( f(Z)\geq t) \geq \mathbb P(|\ell(Z)| \geq t)$ for all norm one linear functionals $\ell$, implies that 
\eqref{eq:Gauss-conc} is sharp 
(up to absolute constants) in the large deviation regime $t>M$ (see also \cite{LMS} and \cite[Proposition 2.9]{PVZ} for details). 
Therefore, in this note the focus is on the one-sided small deviation 
inequality:
\begin{align} \label{eq:Gauss-conc-l}
\mathbb P\left(  f(Z) - M  < -t \right) \leq \frac{1}{2} \exp \left( -\tfrac{1}{2}t^2/L^2 \right), 
\end{align} which holds for all $t>0$ and for any $L$-Lipschitz map $f$. This inequality is of great importance in 
asymptotic geometric analysis, hence one would be interested in refined forms of \eqref{eq:Gauss-conc-l}.
For different ranges of $t$, one can replace the Lipschitz constant $\| \, \|\nabla f\|_2 \, \|_{L_\infty}$ by appropriately chosen 
moments of $\|\nabla f\|_2$.
This is based on various Gaussian functional inequalities such as the logarithmic Sobolev inequalities, the Poincar\'e
inequalities, the $(p,q)$-Poincar\'e inequalities and more (see \cite{BLM,Led}). Even in that case, there exist examples of convex functions 
(e.g. $f(x)=\max_{i\leq n}|x_i|$) for which the $L_2$ norm of the gradient is much larger than the variance, therefore these inequalities fail to 
capture the right order of concentration (see \cite{Cha} for a detailed discussion of this phenomenon). Ideally, one would like to replace the Lipschitz constant in \eqref{eq:Gauss-conc-l} by a statistical measure of dispersion, e.g. the variance. Indeed this is the case for convex functions. Our main result reads as follows: For any convex map $f\in L_2(\gamma_n)$ one has
\begin{align} \label{eq:var-conc}
\mathbb P\left( f(Z)-M <-t \right) \leq \frac{1}{2} \exp \left( -\tfrac{\pi}{1024} t^2/{\rm Var}[f(Z)] \right),
\end{align} for all $t>0$. In view of \eqref{eq:var-ineq} this obviously improves the one-sided concentration inequality in the small
deviation regime. We want to emphasize that the above inequality, unlike to most concentration inequalities which are isoperimetric
in nature, does not follow by the Gaussian isoperimetry. Instead it is obtained by the convexity properties of the Gaussian measure,
thus it could be viewed as a ``new type" of concentration.
The last but not least is that the function is not required to be Lipschitz in \eqref{eq:var-conc}; instead it is valid for 
any convex function $f\in L_2(\gamma_n)$ (in fact we may even prove a similar inequality to \eqref{eq:var-conc} by assuming
weaker integrability condition for $f$; see Remark \ref{rems}.1).

The rest of the paper is organized as follows: In Section 2 we present a proof of the main result. The key ingredient in our argument is
Ehrhard's inequality \cite{Ehr}, inspired by the approach of Kwapien in \cite{Kw}.
We conclude in Section 3 with some applications.

%%%%%%%%%%%%%%%%%%%
\section{Proof of the main result}
%%%%%%%%%%%%%%%%%%%

Let $\Phi$ be the cumulative distribution function of a standard Gaussian random variable, i.e.
\begin{align*}
\Phi(x) =\frac{1}{\sqrt{2\pi}}\int_{-\infty}^x e^{-z^2/2}\, dz, \quad x\in \mathbb R.
\end{align*} Ehrhard's inequality \cite{Ehr} states that for any two convex sets $A,B$ on $\mathbb R^n$ and for any $0<\lambda <1$ one has:
\begin{align*}
\Phi^{-1} [ \gamma_n( (1-\lambda )A +\lambda B)] \geq (1-\lambda )\Phi^{-1}[\gamma_n(A)] + \lambda \Phi^{-1}[\gamma_n(B)].
\end{align*} Ehrhard's result was extended by Lata\l a in \cite{Lat} to the case that one of the two sets is Borel and the other is convex. Finally, in \cite{Bor}, Borell proved that it holds for all pairs of Borel sets. Recently, many different proofs of this fundamental inequality
have appeared in the literature, see e.g. \cite{vHan,IV,NP} and the references therein.

Our goal is to prove the following:

\begin{theorem} \label{thm:main-1}
Let $Z$ be an $n$-dimensional standard Gaussian vector. Let $f$ be a convex function on $\mathbb R^n$ with $f\in L_1(\gamma_n)$ and 
let $M$ be a median for $f(Z)$. Then, we have:
\begin{align*}
\mathbb P \left( f(Z) -M < -t \mathbb E (f(Z)-M)_+ \right) \leq \Phi \left(- \frac{\sqrt{2\pi}}{32} t \right),
\end{align*} for all $t>0$.
\end{theorem}

\noindent {\it Proof.} Since $M$ is a median we have $\mathbb P(f(Z)\leq M)\geq 1/2$. We may assume without loss of generality
that $\mathbb P(f(Z) \leq M)=1/2$. Otherwise we have $\mathbb P(f=M)>0$ and since $f$ is convex we get $f\geq M$, thus
the conclusion is trivially true.
Note that the convexity of $f$ implies that the sub-level sets $\{f \leq t\}, \; t\in \mathbb R$ are convex and the 
function $F(t): = \mathbb P( f(Z) \leq t )$ is log-concave. The latter follows by the following inclusion:
\begin{align*}
(1-\lambda) \{ f \leq t\} +\lambda \{f \leq s\} \subseteq \{ f \leq (1-\lambda)t +\lambda s\},
\end{align*} for $t,s\in \mathbb R$ and $0\leq \lambda \leq 1$ and the fact that $\gamma_n$ is log-concave measure (see \cite[Section 1.8]{Bog} for
the related definition). 
Now, we may use Ehrhard's inequality from \cite{Ehr} (see also \cite[Theorem 4.2.1.]{Bog}) to get that the map $s\mapsto \Phi^{-1} \circ F(s)$, 
$s\in \mathbb R$ is concave (for a proof see \cite[Theorem 4.4.1.]{Bog}). Therefore, we obtain:
\begin{align} \label{eq:ineq-slope}
(\Phi^{-1}\circ F)(M+s) &= (\Phi^{-1} \circ F)(M+s) - (\Phi^{-1} \circ F)(M) \\
&\leq s (\Phi^{-1}\circ F)'(M+) = s\sqrt{2\pi} F'(M+) , \quad s\in \mathbb R. \nonumber
\end{align} Now we give a lower bound for $F'(M+)$ in terms of the standard deviation of $f(Z)$. 

\smallskip

\noindent {\it Claim.} We have the following:
\begin{align*}
F'(M+) \geq \frac{1}{ 32 \mathbb E (f(Z)-M )_+}.
\end{align*}

\smallskip

\noindent {\it Proof of Claim.} Fix $\delta>0$ (that will be chosen appropriately later). Using the log-concavity of $F$ we may write:
\begin{align*}
\delta \frac{F'(M+)}{F(M)} \geq \log F(M+\delta)- \log F(M)  &= \log \big(1 +2 \mathbb P( M < f(Z) \leq M+ \delta) \big) \\
&\geq  \mathbb P \big( M < f(Z) \leq M+ \delta \big) \\
& =  \left( \frac{1}{2} -\mathbb P( f(Z) > M + \delta ) \right),
\end{align*} where we have used the elementary inequality $ \log (1+u )\geq u/2$ for all $0<u\leq 1$. Now we apply Markov's inequality to get:
\begin{align*}
\mathbb P(f(Z) > M +\delta) \leq \frac{\mathbb E (f(Z) -M)_+ }{\delta}.
\end{align*} Combing the above we conclude that:
\begin{align*}
F'(M+) \geq \frac{1}{2\delta} \left( \frac{1}{2} -\frac{ \mathbb E ( f(Z)-M )_+ }{\delta} \right).
\end{align*} The choice $\delta= 4 \mathbb E (f(Z)-M )_+$ yields the assertion of the Claim.

\smallskip 

Going back to \eqref{eq:ineq-slope} we readily see that (for $s=-t \mathbb E(f(Z)-M)_+ $):
\begin{align*}
\Phi^{-1} \left[ \mathbb P \Big( f(Z)-M \leq -t \mathbb E (f(Z)-M)_+ \Big) \right] \leq -t \frac{\sqrt{2\pi} }{ 32},
\end{align*} as required. \prend

\medskip

Let us note that one can prove a similar inequality for the $n$-dimensional exponential measure but for $1$-unconditional functions $f$, i.e.
functions which satisfy $f(x_1,\ldots,x_n)=f(|x_1|,\ldots,|x_n|)$ for all $x=(x_1,\ldots,x_n)\in \mathbb R^n$.  

We fix $W$ for an $n$-dimensional exponential random vector, i.e. 
$W=(\xi_1,\ldots, \xi_n)$, where $(\xi_i)_{i=1}^n$ are independent identically distributed according to the measure $\nu_1$ with density 
function $d\nu_1(x)= \frac{1}{2}e^{-|x|} dx$. Note that if $g_1,g_2$ are i.i.d. standard 
normals and $\xi$ is independent exponential random variable then $|\xi|$ and $\frac{g_1^2+g_2^2}{2}$ have the same distribution (follows easily by checking the moment generating functions). Based on this remark we have the following consequence of Theorem \ref{thm:main-1}: 

\begin{theorem} \label{thm:main-3}
Let $f$ be an 1-unconditional and convex function on $\mathbb R^n$. If $W$ is an exponential random vector on $\mathbb R^n$, then one has:
\begin{align*}
\mathbb P \left( f(W)- M < -t \mathbb E ( f(W)-M)_+  \right)\leq 1-\Phi(ct) \leq \exp(-c' t^2),
\end{align*} for all $t>0$.
\end{theorem}

\noindent {\it Proof.} Consider the function $F:\mathbb R^{2n}\to \mathbb R$ defined as:
\begin{align*}
F(x_1,\ldots,x_n,y_1,\ldots,y_n):=  f\left( \frac{x_1^2+y_1^2}{2}, \ldots, \frac{x_n^2+y_n^2}{2} \right).
\end{align*} Since $f$ is convex and 1-unconditional it follows that $f$ is convex and coordinatewise non-decreasing\footnote{A real valued function $H$ defined on $U \subseteq \mathbb R^k$ is said to be {\it coordinatewise non-decreasing} if it is non-decreasing in each variable while keeping all the 
other variables fixed at any value.} in the octant 
$\mathbb R_+^n= \{z=(z_1,\ldots,z_n) \, : \, z_i\geq 0\}$. Hence $F$ is convex on $\mathbb R^{2n}$. Therefore a direct 
application of Theorem \ref{thm:main-1} yields:
\begin{align*}
\mathbb P \left( f(\tilde W)- M < -t \mathbb E (f(\tilde W)-M)_+ \right) \leq \Phi(-ct),
\end{align*} for all $t>0$, where $\tilde W=(|\xi_1|,\ldots,|\xi_n|)$ and $\xi_i$ are i.i.d. exponential random variables. 
The fact that $f(x_1,\ldots,x_n)=f(|x_1|,\ldots, |x_n|)$ completes the proof. \prend

\begin{remark} \rm In the above argument it is clear that we may also consider longer sums of the form $ g_1^2 +\ldots +g_k^2$. 
That is, if $f:\mathbb R_+^n \to \mathbb R$ is a coordinatewise non-decreasing and convex function, then
\begin{align*}
\mathbb P\left( f(\chi) < M -t \mathbb E(f(\chi)-M)_+ \right)\leq \Phi(-t/2),
\end{align*} for all $t>0$, where $\chi \sim \chi^2(k)$ is a chi squared random variable with $k$ degrees of freedom.
\end{remark}

We conclude this Section with some remarks on the main result.

\begin{remarks} \label{rems} \rm 1. The advantage of this one-sided concentration inequality is that it can be applied for the wide class of 
convex functions which are not necessarily (globally) Lipschitz or which are not even in $L_2(\gamma_n)$; e.g. the function $f(t)=\exp(-t+t^2/2)$ 
is (logarithmically) convex, belongs to $L_1(\gamma_1)$ but $f\notin L_2(\gamma_1)$. Moreover, a careful inspection of the argument shows that 
it is enough to have $f\in L_{1,\infty}(\gamma_n)$ (see e.g. \cite{Graf} for the definition of the weak $L_p$ space) and the conclusion still holds:
\begin{align}
\mathbb P \left( f(Z) < M -t \| (f-M)_+\|_{1,\infty} \right) \leq\Phi(-ct), \quad t>0,
\end{align} where $c>0$ is an absolute constant.\footnote{Here and everywhere else $C$ and $c, c_1, \ldots$ stand for absolute 
constants whose values may change from line to line. We write $c(p)$ if the constant depends only on $p$.}

\smallskip

\noindent 2. Assuming that $\mathbb P(f \leq M)=1/2$, then \eqref{eq:ineq-slope} shows that the variable $f(Z)$ stochastically dominates 
the normal random variable $\zeta:= M + a \cdot g$, where $g$ is a standard normal variable 
and $1/a:= (2\pi)^{1/2} F'(M+) >0$, i.e.
\begin{align*}
\mathbb P( f(Z) \leq s) \leq \mathbb P( \zeta \leq s),
\end{align*} for all $s\in \mathbb R$. Hence one gets $\mathbb E f(Z) \geq \mathbb E\zeta =M$. If $\mathbb P(f\leq M)>1/2$, 
then $\inf f=M$ and the latter is again true. This result is due to Kwapien \cite{Kw}. 
In fact our proof steps on the same starting line as in \cite{Kw}.

\smallskip

\noindent 3. Taking into account the fact that $\mathbb E(f(Z)-M)_+ \leq \mathbb E | f(Z)-M | \leq \sqrt{{\rm Var}f(Z)}$ and 
\begin{align} \label{eq:ineq-phi}
1-\Phi(u)=\Phi(-u)\leq \frac{1}{2}e^{-u^2/2}
\end{align} for all $u>0$ (for a proof see \cite[Lemma 1]{LO}) we immediately get:
\begin{align*}
\mathbb P \left( f(Z)-M < -t \sqrt{ {\rm Var} f(Z)} \right) \leq \Phi \left( -t \frac{\sqrt{2\pi}}{32} \right) \leq \frac{1}{2} \exp\left(- \frac{\pi}{1024} t^2\right), 
\end{align*} for all $t>0$, which is the announced estimate \eqref{eq:var-conc} provided that $f\in L_2(\gamma_n)$.

Furthermore, using the fact $M\geq \mathbb Ef(Z) -\sqrt{{\rm Var}f(Z)}$ once more, we may conclude the following ``Central Limit type'' 
normalization in Theorem \ref{thm:main-1}: 
For any convex function $f$ on $\mathbb R^n$ with $f\in L_2(\gamma_n)$ one has the following distributional inequality:
\begin{align} \label{eq:main-1-clt}
\mathbb P \left( f(Z)-\mathbb Ef(Z) < -t \sqrt{ {\rm Var} f(Z)} \right) \leq \frac{1}{2} \exp \left(- \frac{\pi}{1024} (t-1)^2 \right) < e^{-t^2 /1000},
\end{align} for all $t>1$.

\smallskip

\noindent 4. Let us note that in all the above statements, one can derive the reverse distributional inequality for concave 
functions. Namely, if $f$ is a {\it concave} function on $\mathbb R^n$ with $f\in L_1(\gamma_n)$, then 
\begin{align*}
\mathbb P\left(f(Z)-M >t \mathbb E (M-f(Z))_+ \right)\leq \Phi(-ct),
\end{align*} for all $t>0$, where $M$ is a median for $f(Z)$.

\smallskip

\noindent 5. We should stress the fact that in the statement of Theorem \ref{thm:main-1} we refer to convex functions in $L_1$. 
Thus it is pointless to ask about a similar upper estimate other than the $L_1$-estimate. However in various
significant applications the functions under consideration
are norms or more generally Lipschitz functions which are known to belong in $L_{\psi_2}(\gamma_n)$. 
In fact $\| f-M\|_{\psi_2}\leq C {\rm Lip}(f)$ (where the $L_{\psi_2}$ norm stands for the Orlicz norm with Young function $\psi_2(t)=e^{t^2}-1, \; t\geq 0$). 
However, there are many examples of norms $f$ for which ${\rm Var}[f(Z)] \ll {\rm Lip}(f)^2$. 
Therefore, it is natural to ask if there is one-sided concentration estimate (in the large deviation regime) which takes into 
account both the variance and the Lipschitz constant. A naive approach which puts these remarks together is to combine 
Chebyshev's inequality with the concentration estimate in terms of the Lipschitz constant:
\begin{align*}
\mathbb P( |f(Z)-M | > t ) \leq \exp\left( - \tfrac{1}{2}\max \left\{ \log \left( t/ \sqrt{ {\rm Var} f(Z)} \right), t^2/L^2 \right \} \right).
\end{align*} Even in the case of a norm as above this bound depends continuously on $t>0$ and seems to be the right one.
Example of such a norm is the $\ell_p$ norm on $\mathbb R^n$ with $p=c_0\log n$, for sufficiently small absolute constant $c_0>0$
(see \cite[Section 3]{PVZ}).

\smallskip

\noindent 6. (Non-optimality in $\ell_\infty^n$). Note that Theorem \ref{thm:main-1} for $f(x)=\|x\|_\infty, \; x\in \mathbb R^n$ only yields:
\begin{align*}
\mathbb P( \|Z\|_\infty < (1-\varepsilon)M_{\infty,n} ) \leq \frac{1}{2} e^{-c\varepsilon^2 \log^2 n},
\end{align*} for all $0<\varepsilon <1$, where $M_{\infty,n}$ is the median of $\|Z\|_\infty$. 
This estimate is far from being the sharp one: It is known (see \cite[Claim 3]{Sch2}) that one has:
\begin{align*}
\exp(-C e^{c' \varepsilon \log n}) \leq \mathbb P( \|Z\|_\infty < (1-\varepsilon)M_{\infty,n} ) \leq C\exp(-c e^{c \varepsilon \log n}),
\end{align*} for all $0<\varepsilon<1/2$.

\smallskip

\noindent 7. (Optimality in $\ell_p^n$, $1\leq p<\infty$). In \cite{PVZ} it is proved that for any $1\leq p<\infty$ one has 
$v_{p,n}:={\rm Var}\|Z\|_p/M_{p,n}^2 \leq c(p) /n$, where $M_{p,n}$ is the median for $\|Z\|_p$ 
(see also \cite{PV1} for an extension of this result to any finite dimensional subspace of $L_p$). 
On the other hand, for any norm $\|\cdot\|$ on $\mathbb R^n$ we can deduce that:
\begin{align*}
\mathbb P\left( \|Z\| < (1-\varepsilon) \mathbb E \|Z\| \right) \geq c\exp(-C\varepsilon^2 n),
\end{align*} for all $0<\varepsilon <1/2$. Therefore, we obtain:
\begin{align*}
\mathbb P\left( \|Z\|_p < (1-\varepsilon) M_{p,n} \right)\geq c' \exp \left( -C(p) \varepsilon^2 / v_{p,n} \right).
\end{align*} 

\noindent 8. Probabilistic inequalities similar to \eqref{eq:var-conc}, in the context of log-concave measures, will be presented elsewhere \cite{PV}.

\end{remarks}

%%%%%%%%%%%%%%%%%%%%%%%%%%%%%%%%%%%%%%%
\section{Small ball probabilities and applications}
%%%%%%%%%%%%%%%%%%%%%%%%%%%%%%%%%%%%%%%

In this section we show that the small deviation inequality proved in Theorem \ref{thm:main-1} leads to new reverse H\"older inequalities
for negative moments and small ball probabilities. Toward this end, we exploit once more convexity properties of the Gaussian measure
by utilizing the B-inequality proved by Cordero-Erausquin, Fradelizi and Maurey in \cite{CFM}. The latter states that for any centrally 
symmetric convex body\footnote{A subset $K$ in $\mathbb R^n$ is said to be a centrally symmetric convex body, if it is 
convex, compact with non-empty interior and $K=-K$.} $K$ in $\mathbb R^n$ the function
\begin{align*}
t\mapsto \mathbb P(\|Z\|_K\leq e^t), \quad Z\sim N({\bf 0},I_n)
\end{align*} is log-concave, where $\|\cdot\|_K$ is the gauge of $K$. As this result is available only for norms (the fact that the symmetry 
assumption is essential has been shown in \cite{NT}) from now on we will work within this context. 
Using the aforementioned result, and building on the ideas of Lata\l a and Oleszkiewicz from \cite{LO}, Klartag and Vershynin in \cite{KV} 
introduced a parameter associated with any centrally symmetric convex body which governs the small ball probability for the corresponding
norm. We recall the Klartag-Vershynin parameter (in the Gaussian setting) from \cite{KV}: For any centrally symmetric convex
body $A$ in $\mathbb R^n$ we define
\begin{align*}
d(A):= \min \left \{n, -\log \gamma_n\left(\frac{M}{2}A \right) \right\},
\end{align*} where $M$ is the median of $\|Z\|, \; Z \sim N({\bf 0},I_n)$. Their result reads as follows.

\begin{theorem} [Klartag-Vershynin] \label{thm:KV-sb} Let $A$ be a centrally symmetric convex body in $\mathbb R^n$. Then, one  has
\begin{align*}
 \mathbb P( \|Z\|_A\leq \varepsilon M) \leq 
\frac{1}{2} \varepsilon^{cd(A)}, \quad 0<\varepsilon<1/2,
\end{align*} where $M$ is the median of $\|Z\|_A$ and $Z$
is an $n$-dimensional standard Gaussian vector. 
\end{theorem}

In general it is quite hard to estimate the quantity $d(A)$ and the known lower bounds are in general suboptimal (see Remark \ref{KV-rem}). 
The small deviation inequality from Theorem \ref{thm:main-1} provides a variance-sensitive lower bound for the quantity $d(A)$. 
For this end we associate with any centrally symmetric convex body $A$ in $\mathbb R^n$ the following parameter:
\begin{align}
\beta(A) :=\frac{{\rm Var}\|Z\|_A}{M^2}, \quad  Z\sim N({\bf 0},I_n),
\end{align} where $M$ is the median of $\|Z\|_A$. With this notation we have the following:

\begin{proposition} \label{prop:sb-beta}
Let $A$ be a centrally symmetric convex body in $\mathbb R^n$.
Then, one has the one-sided concentration estimate:
\begin{align} \label{eq:one-side-conc}
\mathbb P\left( \|Z\|_A \leq (1-\varepsilon) M \right) \leq \frac{1}{2} \exp \left(-c \varepsilon^2/\beta(A) \right), \quad 0<\varepsilon<1,
\end{align} where $M$ is the median of $\|Z\|_A$ and 
$Z$ is an $n$-dimensional standard Gaussian random vector. In particular, 
\begin{align} \label{eq: d>beta>k}
d(A) \geq c_1 /\beta(A),
\end{align} therefore, we have the following small ball probability estimate:
\begin{align} \label{eq:s-b-beta}
\mathbb P\left( \|Z\|_A \leq \varepsilon M \right) \leq \frac{1}{2}\varepsilon^{c/\beta(A)},
\end{align} for all $\varepsilon \in (0,1/2)$.
\end{proposition}

\noindent {\it Proof.} We apply Theorem \ref{thm:main-1} for $t=\varepsilon/\sqrt{\beta(A)}$ to get the first estimate.
The bound $d(A) \geq c/\beta(A)$ follows by the definition of $d$ by plugging 
$\varepsilon=1/2$ in \eqref{eq:one-side-conc}. Now the probabilistic estimate \eqref{eq:s-b-beta} follows from 
Theorem \ref{thm:KV-sb} and the obtained lower bound on $d(A)$. \prend

\medskip

It is known that the small ball probability \eqref{eq:s-b-beta} can be easily translated to a small ball probability for 
Gaussian processes (see e.g. \cite[Theorem 7.1]{Led}), thus one has the following formulation.

\begin{theorem}
Let $(G_t)_{t\in T}$ be a centered Gaussian process indexed by a countable set $T$ such that $\sup_{t\in T}|G_t| <\infty$ almost surely. 
Then, for any $\varepsilon\in (0,1/2)$ we have:
\begin{align*}
\mathbb P \left( \sup_{t\in T}|G_t| \leq \varepsilon M \right) \leq \frac{1}{2} \varepsilon^{c M^2/ v^2},
\end{align*} where $M={\rm med}(\sup_{t\in T}|G_t|)$ and $v^2={\rm Var}(\sup_{t\in T}|G_t|)$.
\end{theorem}
 
The proof of the above theorem follows the same lines as in \cite[Theorem 4]{LO} with the obvious adaptions, 
thus it is omitted.
 
\medskip

In view of Theorem \ref{thm:main-3} one can derive small ball estimates for $1$-unconditional norms with respect to the 
exponential measure $\nu_1^n$. This is promised by a result of Cordero-Erausquin, Fradelizi and Maurey, also proved
in \cite{CFM}, that any 1-unconditional log-concave 
measure $\mu$ and 1-unconditional convex body $K$ in $\mathbb R^n$ has the $B$-property, that is $t\mapsto \mu(e^tK)$ 
is log-concave (recently it was proved in \cite{ENT} that the B-property is satisfied by the $\nu_1^n$ and any centrally symmetric convex body). 
Although the proof is the same as in the Gaussian context we sketch it for reader's convenience. 

\begin{proposition} Let $K$ be an 1-unconditional convex body in $\mathbb R^n$. If $W$ is a random vector distributed according 
to the $n$-dimensional exponential measure $\nu_1^n$, then one has 
\begin{align*}
\mathbb P( \|W\|_K \leq \varepsilon m) \leq \frac{1}{2} \varepsilon^{c/\beta}, \quad \varepsilon\in (0,1/2),
\end{align*} where $m$ is the median of 
$\|W\|_K$ and $\beta={\rm Var}\|W\|_K / m^2 $.
\end{proposition}

\noindent {\it Sketch of Proof.} Applying Theorem \ref{thm:main-3} for $x\mapsto \|x\|_K$ we obtain: 
\begin{align} \label{eq:sb-0}
\nu_1^n(\{ x: \|x\|_K \leq m/2 \}) = \nu_1^n \left( \frac{m}{2} K \right) \leq \frac{1}{2}e^{-c/\beta}.
\end{align} On the other hand, since $t\mapsto \nu_1^n(e^tK)$ 
is log-concave, we may argue as follows: given $\varepsilon \in (0,1/2)$ we choose $\lambda \in (0,1)$ 
such that $1/2 = \varepsilon^{1-\lambda}$, i.e.
$1-\lambda=\frac{\log 2}{\log(1/\varepsilon)}$. The log-concavity implies:
\begin{align*}
\nu_1^n \left( \frac{m}{2}K \right) \geq [\nu_1^n( \varepsilon m K)]^{1-\lambda} [\nu_1^n(m K)]^\lambda \quad 
\Longrightarrow
\quad \left[ 2\nu_1^n \left( \varepsilon mK \right)\right]^{1-\lambda} \leq 2\nu_1^n \left( \frac{m}{2}K \right).
\end{align*}
Plug \eqref{eq:sb-0} in the latter we get the assertion. \prend

\medskip

Now we turn in proving reverse H\"older inequalities for negative moments of norms by using
the small deviation \eqref{eq:one-side-conc} and the small ball probability \eqref{eq:s-b-beta}:

\begin{corollary} \label{cor:stab-neg-moms}
Let $K$ be a centrally symmetric convex body in $\mathbb R^n$. Then, one has:
\begin{align*}
\mathbb E \|Z\|_K \left( \mathbb E \|Z\|_K^{-q} \right)^{1/q} \leq \exp \left( C\sqrt{\beta} + C q \beta \right),
\end{align*} for all $0<q< c/\beta(K)$ where $C,c>0$ are absolute constants and $Z$ is an $n$-dimensional standard Gaussian vector.
\end{corollary}

\noindent {\it Proof.} We know that:
\begin{align*}
\mathbb P ( \|Z\|_K \leq \varepsilon M) \leq \frac{1}{2} \varepsilon^{c_1/ \beta} ,  \quad \mathbb P( \|Z\|_K \leq (1-\varepsilon) M) \leq \frac{1}{2}e^{-c_2\varepsilon^2 /\beta}, 
\end{align*} for all $\varepsilon\in (0,1/2)$, where $M$ is the median for $\| Z\|_K$ and $Z\sim N({\bf 0}, I_n)$. Therefore, we may write:
\begin{align*}
\mathbb E\|Z\|_K^{-q} &= M^{-q} \int_0^\infty \mathbb P( \|Z\|_K \leq t M) \frac{q}{t^{q+1}} \, dt \\
& \leq M^{-q} \left( \frac{q}{2} \int_0^{1/2} \varepsilon^{ \frac{c_1}{ \beta } -q-1} \, d\varepsilon + \int_{1/2}^1 \frac{q}{t^{q+1}}P(\|Z\|_K \leq t M)\, dt +1 \right)\\
&\leq M^{-q} \left( \left(\frac{1}{2}\right)^{\frac{c_1}{\beta}-q} \frac{q\beta}{c_1-q\beta} + q \int_0^{1/2} \frac{1}{(1-\varepsilon)^{q+1} } e^{-c_2\varepsilon^2/\beta} \, d\varepsilon +1\right)\\
&\leq M^{-q} \left( 1+ c_3 q \beta + q \int_0^{1/2} \exp( 2(q+1)\varepsilon-c_2\varepsilon^2/\beta) \, d\varepsilon \right),
\end{align*} for all $0<q< c_4/\beta$, where we have also used the elementary inequality $1-u \geq e^{-2u}$ for $0\leq u\leq 1/2$. It is easy to
check that the last integral can be bounded as:
\begin{align*}
\int_0^{1/2} \exp( 2(q+1)\varepsilon-c_2\varepsilon^2/\beta) \, d\varepsilon \leq c_5 \sqrt{\beta} \exp(c_5q^2 \beta),
\end{align*} for all $0< q \leq c_6/\beta$. The result follows. \prend

\begin{remark} \rm \label{KV-rem}
Klartag and Vershynin in \cite{KV} observed that the concentration of measure inequality \eqref{eq:Gauss-conc-l} implies that $d(A)\geq c k(A)$ where $k(A)$  is given by 
\begin{align*}
k(A):= \mathbb E \|Z\|_A^2 /b(A)^2 , \quad b(A)=\max_{\theta \in S^{n-1} } \|\theta\|_A.
\end{align*}
The quantity $k(A)$ is introduced by V. Milman in \cite{Mil} and it is usually referred to as the {\it critical dimension} of the body $A$. 
We refer to \cite{MS} for further information on this quantity. Although the quantity $k(A)$ is easy to be computed, there are several cases in which
bounding $d(A)$ by $k(A)$ gives suboptimal results. Using \eqref{eq:var-ineq} and the fact $\mathbb E\|Z\| \leq c M$ it is clear 
that $ \frac{ 1}{ \beta(A)} \geq c' k(A)$, thus 
Proposition \ref{prop:sb-beta} provides better bounds for the quantity $d(A)$. We illustrate this in the following example: 
Consider as convex body $A$ the unit ball of some $n$-dimensional subspace of $L_{p}, \; 2<p<\infty$.  It is proven in \cite{PV1} that there exists 
a linear image $\tilde A$ of $A$ with $ \beta( \tilde A) \leq C(p) / n$ while $k(\tilde A)$ 
can be of the order $ n^{2/p}$ (up to constants depending only on $p$). In this case the bounds given by Proposition  \ref{prop:sb-beta} are sharp 
(up to constants depending only on $p$). 
\end{remark}

The inequalities presented on the paper can be used to obtain refinements of several classical results in asymptotic 
geometric analysis such as the random version of Dvoretzky's theorem \cite{Mil}. These applications will appear elsewhere \cite{PV}, \cite {PPV}. We close this section by mentioning one interesting application of the results to the Johnson-Linderstrauss flattening lemma.  

The J-L lemma from \cite{JL} (see also \cite{JN}) asserts that: if $\varepsilon \in (0,1)$ and $x_1,\ldots, x_N\in \ell_2$ 
then there exists a linear mapping (which can be chosen to be an orthogonal projection) $P:\ell_2 \to F$, where $F$ is a subspace of $\ell_2$
with $\dim F\leq c\varepsilon^{-2} \log N$ such that
\begin{align*}
(1-\varepsilon)\|x_i-x_j\|_2 \leq \|Px_i-Px_j\|_2 \leq (1+\varepsilon) \|x_i-x_j\|_2,
\end{align*} for all $i,j=1,\ldots, N$.

This dimension reduction principle has found various applications in mathematics and computer science, in addition to the original
application in \cite{JL} for the Lipschitz extension problem. We refer the interested reader to \cite{Ind,KOR,Vem} and the references 
therein for a partial list of its many applications.

The J-L Lemma we are interested in applies for arbitrary target spaces, as was formulated in \cite{Sch3}. Below we suggest a 
refined one-sided version of the latter.

\medskip

\begin{proposition}
Let $X=(\mathbb R^n, \|\cdot\|)$ be a normed space and let $T\subseteq \ell_2^N$ be a finite set with $T=\{u_1,\ldots, u_N\}$. 
The following hold:
\begin{itemize} 

\item [\rm i.] Let $\delta\in (0,1)$ and assume that $\log |T| \lesssim \delta^2 / \beta(X)$. Then, the random Gaussian 
matrix $G=(g_{ij})_{i,j=1}^{n,N}$ satisfies:
\begin{align*}
\|Gu_i- Gu_j\| \geq (1-\delta) \cdot \mathbb E\|Z\| \cdot \|u_i-u_j\|_2,
\end{align*} for all $i,j,=1,\ldots,N$, where $Z\sim N({\bf 0}, I_n)$, with probability greater than $1-ce^{-c \delta^2 / \beta(X)}$.

\item [\rm ii.] Let $\varepsilon\in (0,1/2)$ and assume that  $\log |T| \lesssim \log(1/\varepsilon) /\beta(X)$. 
Then, the random Gaussian matrix $G=(g_{ij})_{i,j=1}^{n,N}$ satisfies:
\begin{align*}
\|Gu_i- Gu_j\| \gtrsim \varepsilon \cdot  \mathbb E\|Z\| \cdot \|u_i-u_j\|_2,
\end{align*} for all $i,j,=1,\ldots,N$, where $Z\sim N({\bf 0}, I_n)$, with probability greater than $1-c\varepsilon^{c/\beta(X)}$.
\end{itemize}
\end{proposition}

\noindent {\it Proof.} Consider $Z_1, \ldots, Z_N$ i.i.d. standard Gaussian vectors on $\mathbb R^n$ and define the 
random matrix $G= [Z_1, \ldots, Z_N]$. Fix $\theta \in S^{N-1}$ and applying 
Theorem \ref{thm:main-1} (as was formulated further in Remark \ref{rems}.3) we get: 
\begin{align*}
\mathbb P( \|G \theta\| < (1- t) \mathbb E\|Z\|) = \mathbb P(\|Z_1\| < \mathbb E\|Z\|-t\mathbb E\|Z\|) \leq C \exp \left( -c t^2 /\beta \right),
\end{align*} for all $t\in (0,1)$. If $T=\{u_1,\ldots, u_N\}$, consider 
the points $\Theta:= \left\{ \frac{u_i-u_j}{\|u_i-u_j\|_2} \, : \, 1\leq i<j\leq N \right\}$ on $S^{N-1}$.
Then, by the union bound we get:
\begin{align*}
\mathbb P( \exists \, \theta \in \Theta \; : \; \|G \theta\| < (1-\delta) \mathbb E\|Z\|) < C_1 N^2\exp(-c_1\delta^2/ \beta) \leq C_2\exp(-c_2 \delta^2 / \beta),
\end{align*} as long as $\log N\leq c \delta^2/\beta(X)$. The assertion follows.

The same reasoning as above, but using \eqref{eq:s-b-beta} instead, yields (ii). \prend

\bigskip

\noindent {\bf Acknowledgments.} The authors are grateful to Ramon van Handel for useful discussions and 
to Mark Rudelson for important remarks. They would also like to thank the anonymous referee whose valuable 
comments helped to improve the presentation of this note.

%%%%%%%%%%%%%%%%%%%%%% references %%%%%%%%%%%%%%%%%%%%%%

%\footnotesize

\medskip

\vspace{.5cm} \noindent 

\begin{minipage}[l]{\linewidth}
  Grigoris Paouris: {\tt grigoris@math.tamu.edu}\\
  Department of Mathematics, Mailstop 3368\\
  Texas A \& M University\\
 College Station, TX 77843-3368\\
  
  \medskip
  
  Petros Valettas: {\tt valettasp@missouri.edu}\\
  Mathematics Department\\
  University of Missouri\\ 
  Columbia, MO 65211\\

\end{minipage}

\end{document}